\documentclass[a4paper,11pt]{article}
\usepackage{kotex}
\usepackage{authblk}
\usepackage{marvosym}
\usepackage{graphicx}
\usepackage{url}
\usepackage{cite}
\usepackage{vmargin}
\setmarginsrb{2cm}{2cm}{2cm}{2cm}{0mm}{0mm}{0mm}{10mm}
\usepackage[colorlinks=true, allcolors=blue, breaklinks=true]{hyperref}
\usepackage{graphicx}
\usepackage{bm}
\usepackage[utf8]{inputenc}
\usepackage[T1]{fontenc}
\usepackage{mathptmx} 
\usepackage{amsmath}
\usepackage{amssymb}
\usepackage{amsthm}
\theoremstyle{definition}
\newtheorem{proposition}{Proposition}
\DeclareMathOperator{\sech}{sech}

\title{\bfseries Harmonic numbers as the summation of integrals}
\author{\normalsize Natanael Karjanto\thanks{\Letter: \url{natanael@skku.edu} \href{https://orcid.org/0000-0002-6859-447X}{\includegraphics[scale=0.08]{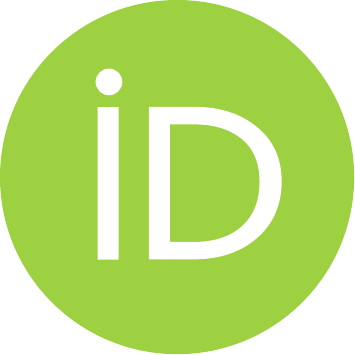}}}}
\affil{Department of Mathematics, University College, Natural Science Campus\\ Sungkyunkwan University, Suwon~16419, Republic of Korea}

\date{\vspace*{-0.5cm} \scriptsize Updated \today}

\begin{document}
\maketitle

\begin{abstract}
Harmonic numbers arise from the truncation of the harmonic series. The $n^\text{th}$ harmonic number is the sum of the reciprocals of each positive integer up to~$n$. In addition to briefly introducing the properties of harmonic numbers, we cover harmonic numbers as the summation of integrals that involve the product of exponential and hyperbolic secant functions. The proof is relatively simple since it only comprises the Principle of Mathematical Induction and integration by parts.
\end{abstract}

\section{Introduction}

The $n^\text{th}$ \emph{harmonic number} is the sum of the reciprocals of the first $n$ natural numbers, expressed in the following form~\cite{graham1994concrete}:
\begin{equation*}
H_n = \sum_{k = 1}^{n} \frac{1}{k} = 1 + \frac{1}{2} + \frac{1}{3} + \frac{1}{4} + \cdots + \frac{1}{n}.
\end{equation*}
The first few harmonic numbers are given as follows:
\begin{align*}
H_1 & = 1 \\
H_2 &= H_1 + \frac{1}{2} = \frac{3}{2} \\
H_3 &= H_2 + \frac{1}{3} = \frac{11}{6} \\
H_4 &= H_3 + \frac{1}{4} = \frac{25}{12} \\
H_5 &= H_4 + \frac{1}{5} = \frac{137}{60} \\
    &\vdots
\end{align*}
Indeed, the harmonic numbers satisfy the following recurrence relation by definition:
\begin{equation*}
H_{n + 1} = H_n + \frac{1}{n + 1}.
\end{equation*}

Since harmonic numbers arise from the truncation of the \emph{harmonic series}, they are the partial sums of it. The divergence of the harmonic series itself can be shown using the integral test by comparing it with the integral of the function $1/x$. Another way to prove its divergence is by comparison test, where each denominator in the series is replaced with the next-largest power of two. A harmonic number can be approximated using the first few terms of the Taylor series expansion~\cite{knuth1997citekey}:
\begin{equation*}
H_n \simeq \gamma + \ln n + \frac{1}{2n},
\end{equation*}
where $\gamma = 0.57721\dots$ is the Euler-Mascheroni constant~\cite{brent1977computation,dence2009survey,havil2017gamma,lagarias2013euler}. See~\cite{owinoh2005effects,caves1996quantum} for applications of this constant. 

The harmonic numbers can be written as Euler's integral representation~\cite{sandifer2007how}:
\begin{equation*}
H_n = \int_{0}^{1} \frac{1 - x^n}{1 - x} \, dx.
\end{equation*}
Harmonic numbers appear in many expressions involving special functions in analytic number theory, including the Riemann zeta function. They are also related to the natural logarithm. For other expositions on harmonic numbers, see~\cite{benjamin2002stirling,choi2011some,spiess1990some,sun2012arithmetic} among others.

In this article, we cover harmonic numbers as the summation of integrals. It provides an extension to one of the proposed problems published in The Problem Corner of a mathematics magazine for students, \emph{The Pentagon}~\cite{karjanto2015problem}, where the solution was published one year later~\cite{karjanto2016solution}. The magazine is the official journal of the Kappa Mu Epsilon (KME), a specialized honor society in mathematics with headquarters at the University of Central Missouri in Warrensburg, Missouri, US~\cite{kappa2021}. KME was founded in 1931 by Emily Kathryn Wyant (1897--1942) at Northeastern Oklahoma State Teachers College. Its focus is to promote the interest in mathematics among undergraduate students. Ever since the first issue appeared in Fall~1941, \emph{The Pentagon} is published semiannually in December and May each year. The purpose of the magazine is not only to serve the needs of college students majoring in mathematics but also to act as a medium through which outstanding student papers on a variety of mathematical topics could be published. In particular, the problems section, known as \emph{The Problem Corner}, contains both proposed and solved problems. The Problem Corner is edited by Pat (Patrick) Costello Kentucky Alpha chapter, Eastern Kentucky University, Richmond, Kentucky, US.

After this introduction, the following section presents the proposition and its proof where harmonic numbers can be expressed as the summation of integrals involving exponential and hyperbolic secant functions. The final section concludes our discussion.

\section{Harmonic numbers as the summation of integrals}		\label{summ}
We have the following proposition.
\begin{proposition}
For $\alpha > 0$ and $n \in \mathbb{N}$, the harmonic number $H_n$ can be represented by the following integral:
\begin{equation*}
H_n = \sum_{k = 1}^n \frac{1}{k} = \frac{1}{2}\sum_{k = 1}^n \left(\int_{-\infty}^\infty e^{-\alpha |x|} \sech^{k + 1} x \, dx + 
\frac{\alpha - (k - 1)}{n} \int_{-\infty}^\infty e^{-(\alpha + 1)|x|} \sech^{k} x \, dx \right) .
\end{equation*}
\end{proposition}

The proof utilizes the Principle of Mathematical Induction (PMI) and integration by parts.
\begin{proof}
Since both integrands of the integrals on the right-hand sides are even functions, then they are equivalent with twice of the integral from $0$ to $+\infty$ (or $-\infty$ to $0$). Without losing generality, we consider the proof by integrating through the former one. Thus, we need to show that
\begin{equation*}
H_n = \sum_{k = 1}^n \frac{1}{k} = \sum_{k = 1}^n \left(I_{k + 1} + \frac{\alpha - (k - 1)}{k} J_{k} \right)
\end{equation*}
where
\begin{align*}
	I_{k + 1} &= \int_{0}^\infty e^{- \alpha x}     \sech^{k + 1} x \, dx,& \alpha > 0, \qquad &k \in \mathbb{N} \\
	J_{k}     &= \int_{0}^\infty e^{-(\alpha + 1)x} \sech^{k} x     \, dx,& \alpha > 0, \qquad &k \in \mathbb{N}.
\end{align*}
We will show it using the PMI. For $k = 1$, we need to show that for $\alpha > 0$
\begin{equation*}
I_2 + \alpha J_1 = 1.
\end{equation*}
Consider $I_2$. By expressing $\sech x = 2 e^x/(1 + e^{2x})$ and letting $y = e^x$, we obtain
\begin{equation*}
I_2 = \int_{0}^\infty e^{-\alpha x} \sech^{2} x \, dx = \int_{1}^\infty \frac{2^2 \, y^{-\alpha + 2}}{(1 + y^2)^2} \frac{dy}{y}
= \int_{1}^\infty 2 y^{-\alpha} \frac{d(1 + y^2)}{(1 + y^2)^2}.
\end{equation*}
Implementing the integration by parts to the last expression and returning to the original variable $x$, it yields
\begin{align*}
I_2 &= \lim_{b \rightarrow \infty} \left.\left(-y^{-(\alpha + 1)} \frac{2y}{1 + y^2} \right)\right|_{1}^{b} - \alpha \int_{1}^\infty y^{-(\alpha + 1)} \frac{2y}{1 + y^2} \, \frac{dy}{y}  \\
	&= \lim_{b \rightarrow \infty} \left. \left(-e^{-(\alpha + 1)x} \sech x \right)\right|_{0}^{b} - \alpha \int_{0}^\infty e^{-(\alpha + 1)x} \sech x \, dx \\
I_2 + \alpha J_1 &= 1, \qquad \alpha > 0.
\end{align*}
Now assume that the inductive step is true for $k = n$, we want to show that it is also true for $k = n +1$.
It is sufficient to show that for $\alpha > 0$ (consequently $\alpha + 1 > 0$)
\begin{equation*}
I_{n + 2} + \frac{\alpha - n}{n + 1} J_{n + 1} = \frac{1}{n + 1}.
\end{equation*}
Employing the similar technique as in the previous case, we observe that
\begin{equation*}
I_{n + 2} = \int_{0}^\infty e^{-\alpha x} \sech^{n + 2} x \, dx 
= \int_{1}^\infty \frac{2^{n + 2} \, y^{-[\alpha - (n + 2)]}}{(1 + y^2)^{n + 2}} \frac{dy}{y}
= \int_{1}^\infty 2^{n + 1} y^{-(\alpha - n)} \frac{d(1 + y^2)}{(1 + y^2)^{n + 2}}.
\end{equation*}
After integrating by parts and returning back to the original variable $x$, we obtain
\begin{align*}
I_{n + 2} &= \lim_{b \rightarrow \infty} \left[-\frac{1}{n + 1} y^{-(\alpha + 1)} \left(\frac{2y}{1 + y^2}\right)^{n + 1} \right]_{1}^{b} - \frac{\alpha - n}{n + 1} \int_{1}^\infty y^{-(\alpha + 1)} \left(\frac{2y}{1 + y^2} \right)^{n + 1} \, \frac{dy}{y}  \\
&= \lim_{b \rightarrow \infty} \left. \left(-\frac{e^{-(\alpha + 1)x}}{n + 1} \sech^{n + 1} x \right)\right|_{0}^{b} - \frac{\alpha - n}{n + 1} \int_{0}^\infty e^{-(\alpha + 1)x} \sech^{n + 1} x \, dx \\
I_{n + 2} + \frac{\alpha - n}{n + 1} J_{n + 1} &= \frac{1}{n + 1}, \qquad \alpha > 0.
\end{align*}
This completes the induction step
\begin{align*}
\textmd{RHS} &= \sum_{k = 1}^{n + 1} \left(I_{k + 1} + \frac{\alpha - (k - 1)}{k} J_{k} \right)  \\
	&= \sum_{k = 1}^n \left(I_{k + 1} + \frac{\alpha - (k - 1)}{k} J_{k} \right) + I_{n + 2} + \frac{\alpha - n}{n + 1} J_{n + 1}\\
	&= \sum_{k = 1}^n \frac{1}{k} + \frac{1}{n + 1} =  \sum_{k = 1}^{n + 1} \frac{1}{k} = \textmd{LHS}.
\end{align*}
Combining with the other part where we integrate from $-\infty$ to $0$, we obtain the desired result. This completes the proof.
\end{proof}

\section{Conclusion}	\label{conclude}

In this article, we have provided a brief coverage of the harmonic numbers. What is particularly interesting is that harmonic numbers can be expressed as the summation of integrals involving the product of exponential and hyperbolic secant functions. By combining the PMI and integration by parts, we have shown that the harmonic numbers can indeed be represented by the summation of the aforementioned integrals.

\section*{Acknowledgments}
The author gratefully acknowledges Dr. Ramadoni for the invitation to contribute this article.	

\subsection*{Dedication}
The author would like to dedicate this article to his late father Zakaria Karjanto (Khouw Kim Soey, 許金瑞) who introduced and taught him the alphabet, numbers, and the calendar in his early childhood. Karjanto senior was born in Tasikmalaya, West Java, Japanese-occupied Dutch~East~Indies on 1~January~1944 (Saturday~Pahing) and died in Bandung, West Java, Indonesia on 18~April~2021 (Sunday~Wage).

\end{document}